# A CHARACTRIZATION OF GRAPHS WITH 3-PATH COVERINGS AND THE EVALUATION OF THE MINIMUM 3-COVERING ENERGY OF A STAR GRAPH WITH M RAYS OF LENGTH 2

PAUL AUGUST WINTER

DEPARTMENT OF MATHEMATICS, UNIVERSITY OF KWAZULU NATAL, HOWARD COLLEDGE, GLENWOOD,DURBAN, SOUTH AFRICA, 4041

email: winter@ukzn.ac.za

ABSTRACT: The smallest set Q of vertices of a graph G, such that every path on 3 vertices, has at least one vertex in Q, is a minimum 3-covering of G. By attaching loops of weight 1 to the vertices of G we can find the eigenvalues associated with G, and hence the minimum 3-covering energy of G. In this paper we characterize graphs with 3-coverings in terms of non-Q-covered edges, and we determine the minimum 3-covering energy of a star graph with m rays each of length 2.

## 1. INTRODUCTION

The Huckel Molecular Orbital theory provided the motivation for the idea of the *energy* of a graph – the sum of the absolute values of the eigenvalues associated with the graph (see [1]). This resulted in the idea of the *minimum 2-covering energy* of a graph in [1] and the minimum 2-covering energy of star graphs with m rays of length 1 were found. This idea was generalized in [2] and the minimum 3-covering energy of complete graphs were determined. In this paper we determine the minimum 3-covering energy of star graphs with m rays of length 2.

## 2. RESULTS OF MINIMUM 3 –COVERINGS OF GRAPHS

All graphs which we shall consider will be finite, simple, loopless and undirected. Let G be such a graph of order n with vertex set $\{v_1, v_2, ..., v_n\}$. A covering (2-covering) of a connected graph G is a set S of vertices of G of such that *every edge* of G has at *least one vertex* in S. Since an edge is a path length 1 on 2 vertices ( a 2-path) we generalize this in [2], in terms of energy, by introducing a *3-covering* (or 3-path covering) of a graph G as being set Q of vertices of G such that every path of G of length 2 (or *3-path*) has *at least one vertex* in Q. Any 3-covering set of G of minimum cardinality is called a *minimum 3-covering* of G.

In the following theorems, Q is a 3-covering of G, and if the vertex x of G, which is not in Q, is a pendant vertex (vertex of degree 1) of G, that belongs to a path P =u,v,w,...,y,x of length s from Q, where u is the only vertex in Q, we say that P is (with respect to Q) a *s-pendant path of G* and the edge yx is with respect to Q) the s-*pendant edge* of P, and y the *middle vertex* of a 3-pendant path of G. If a vertex u is in Q, then the distance of u from Q is taken as 0.

## THEOREM 1

No vertex of G can be a distance of more than 2 from Q.

<u>Proof</u>

Suppose the vertex u of G, that is not in Q, is a distance 3 from Q. Then there exists a path uvwx of length 3 on 4 vertices such that vertices u,v,w are not in Q, but x is in Q. We therefore have a path uvw of length 3 with u,v and w not in Q, which is a contradiction.

**THEOREM 2**

If u is a vertex that is a distance 2 from Q, then u is a pendant vertex.

<u>Proof</u>

Suppose u is on the path uvw of length 2, where vertices u,v are not in Q and w is, and no shorter path exist from u to Q. If u is not a pendant vertex, then it must be connected to a vertex y, where y is not in Q. But then it follows that we have a path vuy on 3 vertices u,v and w which does not have any vertex in Q, a contradiction.

**THEOREM 3**

If uv is an edge of G where u and v do not belong to Q, and neither are pendant vertices, then uv is the middle of a path xuvw P, of length 3 on 4 vertices, such that the ends x and w of the path are both in Q, and/or uv is the edge of the triangle xuv where x is in Q. These edges are disjoint or overlap in both vertices.

<u>Proof</u>

Let uv be and edge of G, neither of which belong to Q or are pendant vertices. Then there must exist vertices w and y such that wuvy is a path on 4 vertices in G. If w is not in Q, then we have a path wuv on 3 vertices with no vertex in Q, a contradiction. Similarly if y is not in Q we get a path uvw on 3 vertices with no vertex in Q. Thus w and y must belong to Q. If w and y are distinct, then we get the 4-path case, otherwise we get the triangle case. Suppose two such edges uv and u'v' have exactly 1 vertex in common, say v'=u'. Then we have a path uv'u' where neither u,v' or u' are in Q, a contradiction. Thus the edges can overlap in both vertices or they must be disjoint.

The path on 4 vertices in theorem 3 is defined as a *3-covering handle-path* of G and the (non-covered )edge uv the *middle edge* of this handle-path.  The triangle in theorem 3 is a *3-covering triangle* of G, and the (non-covered) edge, uv, the *triangle edge* of G. The edge uv in each case is a *non-covered* edge (with respect to Q) of G

## THEOREM 4

If u is a vertex of G then either: (1)  u is in Q ,or (2) if u is not in Q then (2.1) u is a pendant vertex of a 2-pendant path or a 1-pendant path or (2.2) if u is not a pendant vertex, then (2.2.1) u is the middle vertex of a 2-pendant path vuw where v is in Q and w is a pendant vertex  and/or (2.2.2)  u belongs to the path xuv on 3 vertices where x and v are in Q (defined as a *3-covering V path* of G) and/or  (2.2.3) u belongs to the middle edge of a 3-covering handle-path of G and/or u belongs to the triangle edge of a 3-covering triangle of G.

## THEOREM 5

If u is a vertex of G that is not in Q and is a pendant vertex, then u is a pendant vertex of either a 2-pendant path or a 1-pendant path of G. Thus no s-pendant paths exist of length greater than 2.

## THEOREM  6

If uv is an edge of G, where neither u nor v belongs to Q (i.e. a *non-covered edge* of G, then either uv is a pendant edge of a 2-pendant path of G, or uv is the middle edge of a 3-covering handle- path of G, or uv is the triangle edge of a 3-covering triangle of G; the edges must be disjoint, or overlap in both vertices in the case of the non-pendant edges.

Thus there are only **3 types of non-Q-covered edges** of a graph G with a 3-covering set Q- a  2-pendant edge, a middle edge of a handle path, and a triangle

edge – referred to as a 2-pendant, handle or triangle edge. These edges are disjoint except for the edges of the non-pendant kind which can overlap in both vertices. (If a 2-pendant edge uv has v in common with a handle or triangle edge vw, then we will have a path uvw with no vertex in Q.)

**THEOREM 7**

A graph G has a 3-covering Q if and only if the non-Q-covered edges of G are either 2-pendant, handle or triangle edges which are disjoint except for non-pendant edges which can overlap in both vertices.

### 3. MOLECULAR STRUCTURES AND ENERGY

The minimum 2-covering energy of molecular structures given in [1] involves the smallest set of atoms, such that every atom of the structure, is either in the set, or is connected (via bonds) *directly* to at least one vertex of the set. This is generalized to a minimum 3-covering energy of molecular structures, where the smallest set Q of atoms is considered, such that every atom, is either in the set, or connected by a path (of bonded atoms) of *length at most 2*, to at least one atom in the set.

If two atoms u and v are bonded by the edge uv, and neither u and v are in Q, then we say the pair of bonded atoms are *non-Q-covered*. In terms of energy of a structure, in order, say, to prevent destabilization, we may seek the smallest set Q of atoms to be energized, such that *all* non-Q-covered bonded atoms are either, as edges, 2-pendant, handle or triangle edges of the structure (where the edges are disjoint and the non-pendant edges may overlap in both vertices).

### 4. THE MINIMUM-3-COVERING ENERGY OF A GRAPH

A minimum 3-covering matrix of G with a minimum 3-covering set Q of vertices is a matrix:

$$A_Q^3(G) = (a_{i,j})$$

where

$$a_{ij} = \begin{vmatrix} 1 \ \ if \ \ v_i v_j \in E(G) \\ 1 \ \ if \ \ i = j \ \ and \ \ v_i \in Q \\ 0 \ \ otherwise \end{vmatrix} (*)$$

The middle condition (*) is equivalent to loops of weight 1 being attached to the vertices of Q.

The characteristic polynomial of $A_Q^3(G)$ is then denoted by

$$f_n(G, \lambda) \ := \ \det(\lambda I - A_Q^3(G))$$

The *minimum 3-covering energy* (See [1]) is then defined as:

$$E_Q(G) = \sum_1^n |\lambda_i|$$

Where $\lambda_i$ (the *minimum 3-covering eigenvalues*) are the n real roots of the characteristic polynomial.

### 5. THE GENERALIZED STAR GRAPH

The star graph on m+1 verices with m rays of length 1, can be generalized to a star graph with m rays of length n-1:

Take m copies of the path $P_n$, join the paths at their end vertices ,in the centre vertex u,: denote the graph on $n + (n-1)(m-1) = mn - m + 1$ vertices by $S_{1,mP_n}$ ; $m \geq 2, n \geq 2$

If n=2 and $m \geq 3$, then we get the star graph $K_{1,m}$ with m rays of length 1 which has a minimum 3-covering eigenvalues the same as the minimum 2-covering eigenvalues of $K_{1,m}$; the two non-zero eigenvalues are given in [1]:

$$\frac{1+\sqrt{4m+1}}{2}; \frac{1-\sqrt{4m+1}}{2}$$ which implies that the minimum 3-covering energy of $K_{1,m}$ is $\sqrt{4m+1}$.

## 6. THE STAR GRAPH WITH m RAYS OF LENGTH 2

If n=3 and $m \geq 2$, then we label the vertices of the star graph $S_{1,mP_3}$ with m rays of length 2 on 2m+1 vertices as follows:

Centre vertex is u, the set of m vertices a distance 1,2 respectively, from u is labeled:

$$V_1 = \{v_1^1, v_2^1, \ldots v_m^1\}, ; \quad V_2 = \{v_2^1, v_2^2, \ldots v_m^2\}$$ respectively.

The possible 3-covering sets are $\{u\}$; $V_1$- but the minimum 3-covering is the former set.

For constructing the adjacency matrix A of the $S_{1,mP_3}$ we label the center u as $v_1$, the vertices of $V_1, V_2$ as $V_1 = \{v_2, v_3, \ldots v_{m+1}\}, ; \quad V_2 = \{v_{m+2}, v_{m+3}, \ldots, v_{2m+1}\}$ respectively.

For $n = 3, m = 2$ we have the path $P_5$ with minimum 3-covering $Q = \{v_1\}$ with minimum 3-covering adjacency matrix:

$$A_Q^3(S_{1,2P_3}) = \begin{bmatrix} 1 & 1 & 1 & 0 & 0 \\ 1 & 0 & 0 & 1 & 0 \\ 1 & 0 & 0 & 0 & 1 \\ 0 & 1 & 0 & 0 & 0 \\ 0 & 0 & 1 & 0 & 0 \end{bmatrix}$$ so that the characteristic equation is:

$$\det(\lambda I - A_Q^3(S_{1,2P_3})) = \det \begin{bmatrix} \lambda-1 & -1 & -1 & 0 & 0 \\ -1 & \lambda & 0 & -1 & 0 \\ -1 & 0 & \lambda & 0 & -1 \\ 0 & -1 & 0 & \lambda & 0 \\ 0 & 0 & -1 & 0 & \lambda \end{bmatrix}_{5x5}$$ and expanding using

first row:

$$(\lambda-1) \begin{vmatrix} \lambda & 0 & -1 & 0 \\ 0 & \lambda & 0 & -1 \\ -1 & 0 & \lambda & 0 \\ 0 & -1 & 0 & \lambda \end{vmatrix}_{4x4} + \begin{vmatrix} -1 & 0 & -1 & 0 \\ -1 & \lambda & 0 & -1 \\ 0 & 0 & \lambda & 0 \\ 0 & -1 & 0 & \lambda \end{vmatrix} - \begin{vmatrix} -1 & \lambda & -1 & 0 \\ -1 & 0 & 0 & -1 \\ 0 & -1 & \lambda & 0 \\ 0 & 0 & 0 & \lambda \end{vmatrix}_{4x4}$$

Expanding the last 2 matrix determinants about the 3$^{rd}$ and 4$^{th}$ rows:

$$(\lambda-1) \begin{vmatrix} \lambda & 0 & -1 & 0 \\ 0 & \lambda & 0 & -1 \\ -1 & 0 & \lambda & 0 \\ 0 & -1 & 0 & \lambda \end{vmatrix} + \lambda \begin{vmatrix} -1 & 0 & 0 \\ -1 & \lambda & -1 \\ 0 & -1 & \lambda \end{vmatrix}_{3x3} - \lambda \begin{vmatrix} -1 & \lambda & -1 \\ -1 & 0 & 0 \\ 0 & -1 & \lambda \end{vmatrix}_{3x3}$$

Expanding the last two matrix determinants about the 1$^{st}$ and 2$^{nd}$ rows respectively:

$$(\lambda - 1)\begin{vmatrix} \lambda & 0 & -1 & 0 \\ 0 & \lambda & 0 & -1 \\ -1 & 0 & \lambda & 0 \\ 0 & -1 & 0 & \lambda \end{vmatrix} - \lambda\begin{vmatrix} \lambda & -1 \\ -1 & \lambda \end{vmatrix} - \lambda\begin{vmatrix} \lambda & -1 \\ -1 & \lambda \end{vmatrix}$$

The first determinant involves the circulant matrix with solutions:

$$\exp(\frac{2\pi i j}{4})^2 = \exp(\pi i j); \quad n = 0,1,2,3.$$

The second determinant involves the circulant matrix with solutions:

$$\exp(\frac{2\pi i j}{2})^1 = \exp(\pi i j); n = 0,1.$$

Thus the characteristic equation is:

$$(\lambda - 1)^3(\lambda + 1)^2 - 2\lambda(\lambda - 1)(\lambda + 1) = (\lambda - 1)(\lambda + 1)[(\lambda - 1)^2(\lambda + 1) - 2\lambda]$$

$$= (\lambda - 1)(\lambda + 1)[(\lambda^2 - 1)(\lambda - 1) - 2\lambda] = (\lambda - 1)(\lambda + 1)[(\lambda^3 - \lambda^2 - \lambda + 1) - 2\lambda]$$

$$= (\lambda - 1)(\lambda + 1)[(\lambda^3 - \lambda^2 - 3\lambda + 1)]$$

We generalize this to finding the characteristic equation of a star graph on 2m +1 vertices with m rays of length 2:

For $m \geq 2$ and n=3 we have 2m +1 vertices :

$$u, v_1^1, v_2^1, ..., v_m^1, v_1^2, v_2^2, ..., v_m^2$$

$$A(S_{1,mP_3}) = \begin{bmatrix} 1 & 1 & 1 & : & 0 & 0 & 0 \\ 1 & 0 & 0 & 0 & 1 & 0 & 0 \\ : & 0 & 0 & 0 & 0 & 1 & 0 \\ 1 & 0 & 0 & : & 0 & 0 & 1 \\ 0 & 1 & 0 & 0 & 0 & 0 & 0 \\ : & 0 & 1 & 0 & 0 & 0 & 0 \\ 0 & 0 & 0 & 1 & 0 & 0 & 0 \end{bmatrix}_{(2m+1)x(2m+1)}$$
so that the characteristic

equation:

$$\det(\lambda I - A(S_{1,mP_3})) = \det \begin{bmatrix} \lambda-1 & -1 & : & -1 & 0 & : & 0 \\ -1 & \lambda & 0 & 0 & -1 & 0 & 0 \\ : & 0 & : & 0 & 0 & -1 & 0 \\ -1 & 0 & 0 & \lambda & 0 & 0 & -1 \\ 0 & -1 & 0 & 0 & \lambda & 0 & 0 \\ : & 0 & -1 & 0 & 0 & \lambda & 0 \\ 0 & 0 & : & -1 & 0 & 0 & \lambda \end{bmatrix}_{(2m+1)(2m+1)}$$

Expanding the determinant using the first row:

$$= (\lambda-1) \begin{vmatrix} \lambda & 0 & 0 & -1 & 0 & 0 \\ 0 & \lambda & 0 & 0 & -1 & 0 \\ 0 & 0 & \lambda & 0 & 0 & -1 \\ -1 & 0 & 0 & \lambda & 0 & 0 \\ 0 & -1 & 0 & 0 & \lambda & 0 \\ 0 & 0 & -1 & 0 & 0 & \lambda \end{vmatrix}_{(2m)x(2m)}$$

Followed by m determinants:

$$+\begin{vmatrix} -1 & 0 & 0 & -1 & 0 & 0 \\ -1 & \lambda & 0 & 0 & -1 & 0 \\ -1 & 0 & \lambda & 0 & 0 & -1 \\ 0 & 0 & 0 & \lambda & 0 & 0 \\ 0 & -1 & 0 & 0 & \lambda & 0 \\ 0 & 0 & -1 & 0 & 0 & \lambda \end{vmatrix}_{2m x 2m}$$

$$-\begin{vmatrix} -1 & \lambda & 0 & -1 & 0 & 0 \\ -1 & 0 & 0 & 0 & -1 & 0 \\ -1 & 0 & \lambda & 0 & 0 & -1 \\ 0 & -1 & 0 & \lambda & 0 & 0 \\ 0 & 0 & 0 & 0 & \lambda & 0 \\ 0 & 0 & -1 & 0 & 0 & \lambda \end{vmatrix} + ... + \begin{vmatrix} -1 & \lambda & 0 & -1 & 0 & 0 \\ -1 & 0 & \lambda & 0 & -1 & 0 \\ -1 & 0 & 0 & 0 & 0 & -1 \\ 0 & -1 & 0 & \lambda & 0 & 0 \\ 0 & 0 & -1 & 0 & \lambda & 0 \\ 0 & 0 & 0 & 0 & 0 & \lambda \end{vmatrix}$$

Expanding the last m matrix determinants about the $(m+1)^{th}$ , $(m+2)^{th}$ ,...,$2m^{th}$ rows respectively:

$$(\lambda-1)\begin{vmatrix} \lambda & 0 & 0 & -1 & 0 & 0 \\ 0 & \lambda & 0 & 0 & -1 & 0 \\ 0 & 0 & \lambda & 0 & 0 & -1 \\ -1 & 0 & 0 & \lambda & 0 & 0 \\ 0 & -1 & 0 & 0 & \lambda & 0 \\ 0 & 0 & -1 & 0 & 0 & \lambda \end{vmatrix} + \lambda\begin{vmatrix} -1 & 0 & 0 & 0 & 0 \\ -1 & \lambda & 0 & -1 & 0 \\ -1 & 0 & \lambda & 0 & -1 \\ 0 & -1 & 0 & \lambda & 0 \\ 0 & 0 & -1 & 0 & \lambda \end{vmatrix}_{(2m-1)x(2m-1)}$$

$$-\lambda\begin{vmatrix} -1 & \lambda & 0 & -1 & 0 \\ -1 & 0 & 0 & 0 & 0 \\ -1 & 0 & \lambda & 0 & -1 \\ 0 & -1 & 0 & \lambda & 0 \\ 0 & 0 & -1 & 0 & \lambda \end{vmatrix} +..+ \lambda\begin{vmatrix} -1 & \lambda & 0 & -1 & 0 \\ -1 & 0 & \lambda & 0 & -1 \\ -1 & 0 & 0 & 0 & 0 \\ 0 & -1 & 0 & \lambda & 0 \\ 0 & 0 & -1 & 0 & \lambda \end{vmatrix}$$

Last m determinants expanding using row 1,2 …m respectively:

$$(\lambda-1)\begin{vmatrix} \lambda & 0 & 0 & -1 & 0 & 0 \\ 0 & \lambda & 0 & 0 & -1 & 0 \\ 0 & 0 & \lambda & 0 & 0 & -1 \\ -1 & 0 & 0 & \lambda & 0 & 0 \\ 0 & -1 & 0 & 0 & \lambda & 0 \\ 0 & 0 & -1 & 0 & 0 & \lambda \end{vmatrix}_{2mx2m}$$

$$-m\lambda\begin{vmatrix} \lambda & 0 & -1 & 0 \\ 0 & \lambda & 0 & -1 \\ -1 & 0 & \lambda & 0 \\ 0 & -1 & 0 & \lambda \end{vmatrix}_{(2m-2)(2m-2)}$$

The first matrix comes from the circulant matrix with eigenvalues

$$[\exp(\frac{2\pi i j}{2m})]^m; \quad j=0,1,2,..,2m-1, = (-1); m \; times; (1); m \; times$$

The second matrix comes from the circulant matrix with eigenvalues:

$$[\exp(\frac{2\pi ij}{2m-2})]^{m-1}; \quad j = 0,1,2,..,2m-2, = (-1); m-1 \ times; (1); m-1 \ times$$

Which yields the characteristic equation:

$$(\lambda-1)^{m+1}(\lambda+1)^m - m\lambda(\lambda-1)^{m-1}(\lambda+1)^{m-1}$$

$$= (\lambda-1)^{m-1}(\lambda+1)^{m-1}[(\lambda-1)^2(\lambda+1) - m\lambda]$$

$$= (\lambda-1)^{m-1}(\lambda+1)^{m-1}[(\lambda^2-1)(\lambda-1) - m\lambda]$$

$$= (\lambda-1)^{m-1}(\lambda+1)^{m-1}[(\lambda^3 - \lambda^2 - (m+1)\lambda + 1]$$

Thus the graph has minimum 3-covering eigenvalues 1 and -1, each of multiplicity m-1, and 3 eigenvalues from the roots of the cubic equation:

$$\lambda^3 - \lambda^2 - (m+1)\lambda + 1 = 0$$

**THEOREM 8**

The minimum 3-covering energy of the star graph with m rays of length 2 Is:

$$2m - 2 + |x_1| + |x_2| + |x_3| \text{ where } x_1, x_2 x_3 \text{ are the roots of the cubic equation:}$$

$$\lambda^3 - \lambda^2 - (m+1)\lambda + 1 = 0.$$

The 3 roots can be found in (cubic functions in Wikipedia):

$$x_1 = -\frac{b}{3}$$

$$-\frac{1}{3}\sqrt[3]{\frac{1}{2}[2b^3 - 9bc + 27d + \sqrt{(2b^3 - 9bc + 27d)^2 - 4(b^2 - 3c)^3}]}$$

$$-\frac{1}{3}\sqrt[3]{\frac{1}{2}[2b^3-9bc+27d-\sqrt{(2b^3-9bc+27d)^2-4(b^2-3c)^3}\,]}$$

$$x_2=-\frac{b}{3}$$

$$+\frac{1+\sqrt{3}}{6}\sqrt[3]{\frac{1}{2}[2b^3-9bc+27d+\sqrt{(2b^3-9bc+27d)^2-4(b^2-3c)^3}\,]}$$

$$+\frac{1-\sqrt{3}}{6}\sqrt[3]{\frac{1}{2}[2b^3-9bc+27d-\sqrt{(2b^3-9bc+27d)^2-4(b^2-3c)^3}\,]}$$

$$x_3=-\frac{b}{3}$$

$$+\frac{1-\sqrt{3}}{6}\sqrt[3]{\frac{1}{2}[2b^3-9bc+27d+\sqrt{(2b^3-9bc+27d)^2-4(b^2-3c)^3}\,]}$$

$$+\frac{1+\sqrt{3}}{6}\sqrt[3]{\frac{1}{2}[2b^3-9bc+27d-\sqrt{(2b^3-9bc+27d)^2-4(b^2-3c)^3}\,]}$$

With a=1, b=-1, c=-(m+1) and d=1:

$$x_1=\frac{1}{3}$$

$$-\frac{1}{3}\sqrt[3]{\frac{1}{2}[-2-9(m+1)+27+\sqrt{(-2-9(m+1)+27)^2-4(3m+4))^3}\,]}$$

$$-\frac{1}{3}\sqrt[3]{\frac{1}{2}[-2-9(m+1)+27-\sqrt{(-2-9(m+1)+27)^2-4(3m+4)^3}\,]}$$

$$x_2=\frac{1}{3}$$

$$+\frac{1+\sqrt{3}}{6}\sqrt[3]{\frac{1}{2}[-2-9(m+1)+27+\sqrt{(-2-9(m+1)+27)^2-4(3m+4)^3}\,]}$$

$$+\frac{1-\sqrt{3}}{6}\sqrt[3]{\frac{1}{2}[-2-9(m+1)+27-\sqrt{(-2-9(m+1)+27)^2-4(3m+4)^3}\,]}$$

$$x_3 = \frac{1}{3}$$

$$+ \frac{1-\sqrt{3}}{6} \sqrt[3]{\frac{1}{2}[-2-9(m+1)+27+\sqrt{(-2-9(m+1)+27)^2 - 4(3m+4)^3}\,]}$$

$$+ \frac{1+\sqrt{3}}{6} \sqrt[3]{\frac{1}{2}[-2-9(m+1)+27-\sqrt{(-2-9(m+1)+27)^2 - 4(3m+4)^3}\,]}$$

Simplifying:

$$x_1 = \frac{1}{3}$$

$$- \frac{1}{3} \sqrt[3]{\frac{1}{2}[16-9m+\sqrt{27m^3+189m^2-144m+308}\,]}$$

$$- \frac{1}{3} \sqrt[3]{\frac{1}{2}[16-9m-\sqrt{27m^3+189m^2-144m+308}\,]}$$

$$x_2 = \frac{1}{3}$$

$$+ \frac{1+\sqrt{3}}{6} \sqrt[3]{\frac{1}{2}[16-9m+\sqrt{27m^3+189m^2-144m+308}\,]}$$

$$+ \frac{1-\sqrt{3}}{6} \sqrt[3]{\frac{1}{2}[16-9m-\sqrt{27m^3+189m^2-144m+308}\,]}$$

$$x_3 = \frac{1}{3}$$

$$+ \frac{1-\sqrt{3}}{6} \sqrt[3]{\frac{1}{2}[16-9m+\sqrt{27m^3+189m^2-144m+308}\,]}$$

$$+\frac{1+\sqrt{3}}{6}\sqrt[3]{\frac{1}{2}[16-9m-\sqrt{27m^3+189m^2-144m+308}\,]}$$